\documentclass{amsart}
\usepackage{amssymb}
\usepackage[all]{xy}
\newtheorem{theorem}{Theorem}[section]
\newtheorem{lemma}[theorem]{Lemma}

\newtheorem{cor}[theorem]{Corollary}

\theoremstyle{definition}

\theoremstyle{remark}

\numberwithin{equation}{subsection}

\def\qed{$\Box$}

\def\ff{\mathbb{F}}

\def\qq{\mathbb{Q}}

\def\Tr{\mathrm{Tr}}
\def\Nm{\mathrm{Norm}}

\begin{document}

\title{On Katz's bound for number of elements with given trace and norm}

\author{Marko Moisio}
\address{Department of Mathematics and Statistics, University of
Vaasa, P.O. Box 700, FIN-65101, Vaasa, Finland}
\address{ Email: mamo@uwasa.fi}

\author{Daqing Wan}
\address{Department of Mathematics, University of California,
Irvine, Ca92697-3875}
\address{Email: dwan@math.uci.edu}
\date{}



\begin{abstract}
In this note an improvement of the Katz's bound on the number of elements in a finite field with given trace and norm is given. The improvement is obtained by reducing the problem to estimating the number of rational points on certain toric Calabi-Yau hypersurface, and then to use detailed cohomological calculations by Rojas-Leon and the second author for such toric hypersurfaces.
\end{abstract}

\maketitle
\setcounter{tocdepth}{1}
\section{Introduction}
Let $p$ be a prime and $\ff_q$ be the finite field of $q$ elements
of characteristic $p$. Given $a, b\in \ff_q$, and positive integer
$m\geq 2$, let
$$N_m(a,b)=\#\{ \alpha\in \ff_{q^m} |
{\rm Tr}_{\ff_{q^m}/\ff_q}(\alpha)=a, {\rm
Norm}_{\ff_{q^m}/\ff_q}(\alpha)=b\}.$$ Motivated by various
applications, it is of interest to give a sharp estimate for the
number $N_m(a,b)$. The case $b=0$ is trivial.

Katz \cite{Ka} proved the following bound:

\begin{theorem}Let $a, b \in {\ff^*_q}$ and $n\geq 1$. Then
$$|N_{n+1}(a,b) -\frac{q^{n+1}-1}{q(q-1)}| \leq
(n+1)q^{\frac{n-1}{2}}.$$
\end{theorem}
This bound was used by Moisio \cite{Mo} to improve some cases of the
explicit bound in Wan \cite{Wa2} on the number of irreducible
polynomials in an arithmetic progression of $\ff_q[x]$. In the case
$n+1=3$, the Katz bound also plays a significant role in Cohen and
Huczynska \cite{Ch} for their proof of the existence of a cubic
primitive normal polynomial with given norm and trace.

If $a=0$, Katz's bound can be improved in an elementary way using
character sums \cite{Mo}:
$$|N_{n+1}(0,b) -\frac{q^{n}-1}{q-1}| \leq
(d-1)q^{\frac{n-1}{2}},$$ where $d=\gcd(n+1, q-1)$.

In this note, we give a uniform improvement of Katz's bound in the
case $a\not=0$.

\begin{theorem}Let $a, b \in {\ff^*_q}$ and $n\geq 1$. Then
$$|N_{n+1}(a,b) -\frac{q^{n}-1}{q-1}| \leq
nq^{\frac{n-1}{2}}.$$
\end{theorem}

In the case that $n+1$ is a power of $p$, this improvement was first
proved by Moisio \cite{Mo} using Deligne's estimate for
hyper-Kloosterman sums. Moreover, in the case $n+1=3$ also the bounds
\[
    3\left\lceil\frac{q+1-2\sqrt q}3\right\rceil\le N_3(a,b)
    \le3\left\lfloor\frac{q+1+2\sqrt q}3\right\rfloor.
\]
were obtained in \cite{Mo} by using the Hasse's bound for elliptic
curves together with a divisibility result. In corollary 2.4, we
extend such divisibility bounds to $N_{\ell}(a,b)$, where $\ell\geq
3$ is any prime.

In the general case, our proof of Theorem 1.2 consists of two steps.
The first step is to reduce it to estimating the number of
$\ff_q$-rational points on certain toric Calabi-Yau hypersurface
over $\ff_q$. The second step is to use the detailed cohomological
calculations in Rojas-Leon and Wan \cite{RW07} for such toric
hypersurfaces. In the case $n+1=3$, the above improved bounds should
significantly reduce the amount of calculations in \cite{Ch}.


\section{Proof of Theorem 1.2}

Let $u=b/a^{n+1} \in \ff^*_q$. Let $N(u)$ denote the number of
$\ff_q$-rational points on the toric hypersurface
$$Y_u: \ X_1 +\cdots + X_n +\frac{u}{X_1\cdots X_n} -1 =0.$$

\begin{lemma}
$$N_{n+1}(a,b) = \frac{q^{n}-1}{q-1}
+(-1)^n\left( N(u)-\frac{(q-1)^n-(-1)^n}{q}\right).$$
\end{lemma}

{\sl Proof}. Write the equation of $Y_u$ in the form
\begin{eqnarray*}
    &X_1+\cdots+X_{n+1}&=1\\
    &X_1\cdots X_{n+1}&=u.
\end{eqnarray*}
Let $\psi$ be the canonical additive character of $\ff_q$. Now
\[
q(q-1)N(u)=\sum_{x_1,\dots,x_{n+1}}\sum_v\psi(v(x_1+\cdots+x_{n+1}-1))
\sum_{\chi}\chi(u^{-1}x_1\cdots x_{n+1}),
\]
where $x_1,\dots,x_{n+1}$ run over $\ff_q^*$, $v$ runs over $\ff_q$,
and $\chi$ runs over the multiplicative character group of $\ff_q$.

Let $G(\chi)$ denote the Gauss sum
\[
    G(\chi)=\sum_{x\in\ff_q^*}\psi(x)\chi(x).
\]
It follows that
\begin{eqnarray}\label{e:N_u}
q(q-1)N(u)
&=&(q-1)^{n+1}+\sum_{v\ne0}\psi(-v)\sum_{\chi}\bar\chi(u)
\prod_{i=1}^{n+1}\sum_{x_i}\psi(vx_i)\chi(x_i)\notag\\
&\stackrel{x_i\mapsto {x_i/v}}{=}&(q-1)^{n+1}+\sum_{v\ne0}\psi(-v)
\sum_{\chi}\bar\chi(uv^{n+1})G(\chi)^{n+1}\notag\\
&=&(q-1)^{n+1}+\sum_{\chi}G(\chi)^{n+1}\bar\chi(u)
\sum_{v\ne0}\psi(-v)\bar\chi^{n+1}(v)\notag\\
&=&(q-1)^{n+1}+\sum_{\chi}G(\chi)^{n+1}G(\bar\chi^{n+1})\bar\chi((-1)^{n+1}u).
\end{eqnarray}

Next we express $N_{n+1}(a,b)$ in terms of Gauss sums. We use the
abbreviated notations $\Tr$ and $\Nm$ in place of
$\Tr_{\ff_{q^{n+1}}/\ff_q}$ and $\Nm_{\ff_{q^{n+1}}/\ff_q}$. Let
$\psi_{n+1}=\psi\circ {\Tr}$ be the canonical additive character of
$\ff_{q^{n+1}}$ and let $\alpha\in\ff_{q^{n+1}}$ with
$\Tr(\alpha)=1$. Now,
\begin{eqnarray*}
q(q-1)N_{n+1}(a,b)&=&\sum_{x\in\ff_{q^{n+1}}^*}\sum_v\psi(v(\Tr(x-\alpha
a))
\sum_{\chi}\chi(b^{-1}\Nm(x))\\
&=&\sum_v\psi(-av)\sum_{\chi}\bar\chi(b)\sum_x \psi_{n+1}(vx)\chi(\Nm(x))\\
&=&q^{n+1}-1+\sum_{v\ne0}\psi(-av)\sum_{\chi}\bar\chi(b)\sum_x \psi_{n+1}(vx)\chi(\Nm(x))\\
&\stackrel{x\mapsto
x/v}{=}&q^{n+1}-1+\sum_{v\ne0}\psi(-av)\sum_{\chi}\bar\chi(bv^{n+1})
\sum_x \psi_{n+1}(x)\chi(\Nm(x)),
\end{eqnarray*}
since $\Nm(v)=v^{n+1}$.

By the Davenport-Hasse identity the inner sum
\[
    \sum_x \psi_{n+1}(x)\chi(\Nm(x))=(-1)^nG(\chi)^{n+1},
\]
and therefore
\begin{eqnarray*}
q(q-1)N_{n+1}(a,b)&=&q^{n+1}-1+(-1)^n\sum_{\chi}G(\chi)^{n+1}\bar\chi(b)
\sum_{v\ne0}\psi(-av)\bar\chi^{n+1}(v)\\
&=&q^{n+1}-1+(-1)^n\sum_{\chi}G(\chi)^{n+1}G(\bar\chi^{n+1})\bar\chi((-1)^{n+1}b/a^{n+1}).
\end{eqnarray*}
Comparing this expression with~(\ref{e:N_u}), one finds that
$$N_{n+1}(a,b) = \frac{q^{n+1}-1}{q(q-1)}
+(-1)^n\left( N(u)-\frac{(q-1)^{n+1}}{q(q-1)}\right).$$ One checks
that this is the same as the expression in Lemma 2.1. \hfill\qed

This lemma reduces Theorem 1.2 to the following

\begin{theorem} Let $u\in \ff^*_q$. Then
$$|N(u)-\frac{(q-1)^n-(-1)^n}{q}| \leq n q^{\frac{n-1}{2}}.$$
\end{theorem}

{\sl Proof}. Over the algebraic closure $\bar{\ff}_q$, we can write
$u=\lambda^{-(n+1)}$ for some non-zero element $\lambda$. Then $Y_u$
is isomorphic to the toric hypersurface
$$X_{\lambda}: \ X_1 +\cdots +X_n + \frac{1}{X_1\cdots X_n} -\lambda
=0$$ whose zeta function over a finite field was studied in detail
in \cite{RW07}, see \cite{Wa} for more elementary description of the
results. For a prime $\ell\not=p$, the $\ell$-adic cohomology
$$H_c^j(Y_u\otimes \bar{\ff}_q, \qq_{\ell})\cong H_c^j(X_{\lambda}\otimes \bar{\ff}_q, \qq_{\ell})
$$
was calculated in Theorem 2.1 in \cite{RW07}. In particular, we have
$$H_c^j(Y_u\otimes \bar{\ff}_q, \qq_{\ell}) =0, \ j<n-1 \ {\rm or} \
j>2n-1,$$
$$H_c^j(Y_u\otimes \bar{\ff}_q, \qq_{\ell}) \cong \qq_{\ell}^{{n\choose
j-n+2}}(n-1-j), \ n\leq j\leq 2n-2,$$ and there is an exact sequence
of Galois modules
$$0\rightarrow \qq_{\ell}^n \rightarrow H_c^{n-1}(Y_u\otimes \bar{\ff}_q,
\qq_{\ell})\rightarrow M_u\rightarrow 0,$$ where $M_u$ is of rank at
most $n$ and mixed of weight at most $n-1$. It follows that
$$|{\rm Tr}({\rm Frob}_u|M_u)| \leq nq^{\frac{n-1}{2}}.$$
By the $\ell$-adic trace formula,
$$N(u) =\sum_{j=n}^{2n-2}(-1)^j{n\choose j-n+2} q^{(j-(n-1))}
+(-1)^{n-1}n+(-1)^{n-1}{\rm Tr}({\rm Frob}_u|M_u).$$ Replacing $j$
by $j+n-2$, one finds
$$N(u) =\sum_{j=2}^n (-1)^{j-n}{n\choose j}q^{(j-1)} +(-1)^{n-1}n+(-1)^{n-1}{\rm
Tr}({\rm Frob}_u|M_u).$$ The theorem follows. \hfill\qed

{\bf Remark}. If $u\not= (n+1)^{-(n+1)}$,  i.e., $\lambda \not\in \{
(n+1)\zeta |\zeta^{n+1}=1\}$, then $M_u$ is pure of weight $n-1$ and
of rank $n$. If $u= (n+1)^{-(n+1)}$ (necessarily $p\not| n+1$), then
the rank of $M_u$ drops by $1$ and thus
$$|{\rm Tr}({\rm Frob}_u|M_u)| \leq (n-1)q^{\frac{n-1}{2}}.$$
If $u= (n+1)^{-(n+1)}$ and $n$ is even, then one of the Frobenius
eigenvalues has weight $n-2$ (instead of $n-1$), and thus
$$|{\rm Tr}({\rm Frob}_u|M_u)| \leq (n-2)q^{\frac{n-1}{2}}+q^{\frac{n-2}{2}}
.$$ All these follow from Proposition 2.6 in \cite{RW07}.

\begin{cor}Let $u= (n+1)^{-(n+1)}$. Then
$$|N(u)-\frac{(q-1)^n-(-1)^n}{q}| \leq (n-1) q^{\frac{n-1}{2}}.$$
If $n$ is also even, then
$$|N(u)-\frac{(q-1)^n-(-1)^n}{q}| \leq (n-2) q^{\frac{n-1}{2}}+q^{\frac{n-2}{2}}.$$
\end{cor}

\begin{cor}Let $\ell\geq 3$ be a prime number. Let $a,b\in \ff^*_q$.
Then, we have
$$\ell\left\lceil\frac{\frac{q^{\ell -1}-1}{q-1}-(\ell-1)q^{(\ell-2)/2}}\ell\right\rceil\le N_{\ell}(a,b)
    \le \ell\left\lfloor\frac{\frac{q^{\ell -1}-1}{q-1}+(\ell-1)q^{(\ell-2)/2}}\ell\right\rfloor.$$
\end{cor}

{\sl Proof}. Let $R$ be the number of $c\in \ff_q$ such that $\ell c
=a$ and $c^{\ell}=b$. It is clear that $R$ is either $0$ or $1$.
Since $\ell$ is a prime, $N_{\ell}(a,b)-R$ is divisible by $\ell$.
If $R=0$, the corollary is the consequence of Theorem 1.2 and the
divisibility  of $N_{\ell}(a,b)$ by $\ell$.

Assume now that $R=1$. Since $a\not=0$, $\ell$ cannot be $p$. In
this case, we have $a=\ell c$, $b=c^{\ell}$ and thus $u=b/a^{\ell} =
\ell^{-\ell}\in \ff^*_q$. We can apply the stronger estimate in the
previous corollary to deduce the desired inequalities for
$N_{\ell}(a,b)$.

\end{document}